\documentclass[10pt,a4paper]{article}
\usepackage[utf8]{inputenc}
\usepackage{amsmath, amssymb} 
\newcommand{\K}{\operatornamewithlimits{K}}

\usepackage{ccfonts} 

\begin{document}
\title{Some conjectural continued fractions}
\author{Thomas \textsc{Baruchel}\thanks{baruchel@riseup.net}}
\date{\today}

\maketitle
\section*{Method and notation}
All results presented below have been found by experimental means. Millions of random continued fractions have been computed at very high speed (code was written in~C and double-precision floating numbers were used for the greatest part of the computation despite their relatively low precision). Relations were detected by matching values against various databases with a fast implementation of the famous PSLQ algorithm. Of course, all expressions have been later checked with a much higher precision, but they are not proved and they are copied here as mere conjectures.

Several databases of interesting constants were used; one of them was built from the book of Michael Shamos. Later {\em ad hoc} databases were built in order to focus more accurately on specifical constants.

For convenience reasons, the notation being used below for the generalized continued fractions is the old one used by Gauss: $\K_{n=1}^\infty \frac{b_n}{a_n}$. It has to be remembered however that the expression at the right of the~$\K$ operator is not a complete fraction that could be cancelled down.
\[
a_0+\displaystyle\K_{n=1}^\infty\displaystyle\frac{b_n}{a_n}\;=\;
a_0+
  \genfrac{}{}{0pt}{0}{}{}
  \rlap{$\dfrac{b_1}{\phantom{a_1}}$}%
  \genfrac{}{}{0pt}{0}{}{a_1+
    \genfrac{}{}{0pt}{0}{}{}
    \rlap{$\dfrac{b_2}{\phantom{a_2}}$}%
    \genfrac{}{}{0pt}{0}{}{a_2+
      \genfrac{}{}{0pt}{0}{}{}
      \rlap{$\dfrac{b_3}{\phantom{a_3}}$}%
      \genfrac{}{}{0pt}{0}{}{a_3+
        \genfrac{}{}{0pt}{0}{}{\ddots}
  }}}%
\]

\section{A continued fraction involving $\Gamma\left(\xi+1/2\right)/\Gamma\left(\xi\right)$}
A continued fraction for the quotient~$\Gamma\left(\xi+1/2\right)/\Gamma\left(\xi\right)$ involving simple polynomial functions was conjectured to exist and therefore an algorithm was written for gathering as many constants related to the~$\Gamma$~function as possible. Many identities were detected for particular values and a general expression was found:
\begin{equation}\label{gamma1}
    \K_{n=1}^\infty\displaystyle\frac
  {\left(\alpha n+1\right)\left(-2\alpha n+2\alpha\xi-1\right)}
  {3\alpha n+2+\alpha}
    \;=\;
   -1-\alpha +
\displaystyle\frac
{\Gamma\left(\xi+\displaystyle\frac{1}{2}\right)\Gamma\left(\displaystyle\frac{1}{2\alpha}\right)}
{\Gamma\left(\xi\right)\Gamma\left(\displaystyle\frac{\alpha+1}{2\alpha}\right)}
\end{equation}

\section{An identity involving the lemniscate constant}

By looking for values related to the lemniscate constant in randomly generated continued fractions made of integer terms, a general expression could be found:
\begin{equation}\label{gamma2}
    \begin{array}{l}
        \displaystyle
    \K_{n=1}^\infty\displaystyle\frac
    {-\left(2n+\alpha\right)\left(4n-4\xi+2\alpha-1\right)}
    {6n+2\alpha+2}
    \\[16pt]\qquad =\;
   -2-\alpha +
\displaystyle\frac
        { 2\sqrt{2}\;\Gamma\left(\displaystyle\frac{1}{4}\right)
        \Gamma\left(\displaystyle\frac{3}{4}\right)
      \Gamma\left(1+\displaystyle\frac{\alpha}{4}\right)
      \Gamma\left(\displaystyle\frac{3}{4}-\frac{\alpha}{4}+\xi\right) }
    { \pi\;
      \Gamma\left(\displaystyle\frac{1}{2}+\frac{\alpha}{4}\right)
      \Gamma\left(\displaystyle\frac{1}{4}-\frac{\alpha}{4}+\xi\right) }
    \end{array}
\end{equation}
where the case~$\alpha=1$ here covers the case~$\alpha=2$ in the previous formula~\eqref{gamma1}.


\section{A functional relation involving self powers}
The function~$g$ is defined as a generalized continued fraction
\[
    g\left(k,x\right) = 
\displaystyle\operatorname*{K}_{n=1}^{\infty}
    \frac
    {\left(n+1\right)\left(k-n-k/x\right)/x}
    {\left(n+1\right)\left(1+1/x\right)}
    \]
then the following relation stands:
\begin{equation}\label{funcrel}
    \begin{array}{l}
        2\;+\;
    g\left(\alpha,\xi\right) \;+\; g\left(\alpha,\displaystyle\frac{\xi}{\xi-1}\right)
        \\[12pt]
        \qquad\qquad =\;
        \alpha\left(\xi-1\right)^{\alpha/\xi-1}\left(\displaystyle\frac{\xi}{\xi-1}\right)^{\alpha-2} \textrm{B}\left(\alpha/\xi, \alpha-\alpha/\xi\right)
    \end{array}
\end{equation}
where~$\textrm{B}$ is the beta function.

An easy case occurs when~$\alpha=\xi/\left(\xi-1\right)$, which leads to~$g\left(\alpha,\xi\right)=0$. 
Then, choosing~$\xi=\phi$ (the golden ratio) gives the nice identity:
\begin{equation}
    \phi^\phi
    \;=\;
    2 +
\displaystyle\operatorname*{K}_{n=1}^{\infty}
    \frac
    {\left(n+1\right)\left(1-n/\phi\right)/\phi}
    {\left(n+1\right)\left(2-1/\phi\right)}
\end{equation}    
while choosing~$\xi=\phi^2$ gives:
\begin{equation}
  \phi^{2/\phi}
    \;=\;
    2 +
\displaystyle\operatorname*{K}_{n=1}^{\infty}
    \frac
    {\left(n+1\right)\left(1-\left(n+1\right)/\phi\right)}
    {\left(n+1\right)\phi}
\end{equation}
and similarly for other constants like~$\sqrt{3}$, $\sqrt[3]{4}$, $\sqrt[4]{5}$, etc.
The general expression coming from this case~$\alpha=\xi/\left(\xi-1\right)$ being:
\begin{equation}\label{selfpower}
    \left(\displaystyle\frac{x}{x-1}\right)^{x-1}
    \;=\; 2 +
\displaystyle\operatorname*{K}_{n=1}^{\infty}
    \frac
    {\left(n+1\right)\left(x-n-1\right)/x}
    {\left(n+1\right)\left(x+1\right)/x}
\end{equation}

Building a similar identity relying on~$g\left(2\xi,\xi/\left(\xi-1\right)\right)$ is not very difficult, but the resulting expression is not as simple as previously:
\begin{equation}
    \displaystyle\frac{x-1}{2x-1}\left(
      \left(\displaystyle\frac{x}{x-1}\right)^{2x-1} - 1
    \right)
    \;=\;2+
\displaystyle\operatorname*{K}_{n=1}^{\infty}
    \frac
    {\left(n+1\right)\left(2x-n-2\right)/x}
    {\left(n+1\right)\left(x+1\right)/x}
\end{equation}

\section{A continued fraction for $\tanh z$}
While trying to detect miscellaneous constants by building continued fractions with simple polynomials, several different values related to the hyperbolic tangent function were detected. Generalizing happened to be very easy and the resulting continued fraction seems to be different from the usual one:
\begin{equation}
\tanh z
  \;=\;
  \K_{n=1}^\infty\displaystyle\frac
    {1 + \left(\left(n-1\right)\displaystyle\frac{\pi}{4}z^{-1}\right)^2_{\vphantom{\frac{X}{X}}}  }
    {\left(2n-1\right)\displaystyle\frac{\pi^{\vphantom{a}}}{4}z^{-1} }
\end{equation}

Convergence is slower however than the classical continued fraction for $\tanh$. On the other hand, this version produces nice and simple continued fractions when used with multiples of~$\pi$.
%
%
%

\section{A continued fraction for a sum of products}
Though Euler's continued fraction formula gives an easy continued fraction for the following sum of products, the continued fraction below seems to have a quicker convergence:
\begin{equation}
\displaystyle\sum_{n=1}^\infty\displaystyle\frac{1}
  {\prod_{k=0}^n \alpha k + 1}
  \;=\; \displaystyle\K_{n=1}^\infty \displaystyle\frac{\left(\alpha n\right)^{-1}}{1}
\end{equation}

Interesting values occur for several values of $\alpha$:
{\small
\[
\begin{array}{cl@{\qquad\qquad}cl}
-\displaystyle\sqrt{\displaystyle\frac{\pi}{2e}}\;\textrm{erfi}\,\frac{1}{\sqrt{2}}
  &\textrm{for } \alpha=-2
  &
\displaystyle\frac{e^2-5}{2}
  &\textrm{for } \alpha=\displaystyle\frac{1}{2}\\[16pt]
e-2\vphantom{\displaystyle\frac{e^2-5}{2}}
  &\textrm{for } \alpha=\displaystyle 1
  &
-1+\displaystyle\sqrt{\displaystyle\frac{e\pi}{2}}\;\textrm{erf}\,\frac{1}{\sqrt{2}}
  &\textrm{for } \alpha=2\\[16pt]
\end{array}
\]
} 

\end{document}